# ISS-Based Robustness to Various Neglected Damping Mechanisms for the 1-D Wave PDE

**Iasson Karafyllis[*] and Miroslav Krstic[**]**

[*]Dept. of Mathematics, National Technical University of Athens, Zografou Campus, 15780, Athens, Greece, email: iasonkar@central.ntua.gr

[**]Dept. of Mechanical and Aerospace Eng., University of California, San Diego, La Jolla, CA 92093-0411, U.S.A., email: krstic@ucsd.edu

**Abstract**

This paper is devoted to the study of the robustness properties of the 1-D wave equation for an elastic vibrating string under four different damping mechanisms that are usually neglected in the study of the wave equation: (i) friction with the surrounding medium of the string (or viscous damping), (ii) thermoelastic phenomena (or thermal damping), (iii) internal friction of the string (or Kelvin-Voigt damping), and (iv) friction at the free end of the string (the so-called passive damper). The passive damper is also the simplest boundary feedback law that guarantees exponential stability for the string. We study robustness with respect to distributed inputs and boundary disturbances in the context of Input-to-State Stability (ISS). By constructing appropriate ISS Lyapunov functionals, we prove the ISS property expressed in various spatial norms.

## 1. Introduction

The wave equation is the prototype Partial Differential Equation (PDE) for the description of vibrations in elastic solid and compressible fluid media. Among the numerous physical phenomena that it models are lateral, longitudinal, and torsional oscillations in strings and acoustic oscillations (in gases, liquids and solids; see for instance the textbooks [13,29]). The study of the dynamics of the wave equation with or without boundary control has attracted the interest of many researchers (see for instance [3,4,5,6,7,11,12,14,15,17,19,20,21,22,24,25,26,27,30,31]). However, the wave equation is an "idealized" model. In practice there are various damping mechanisms that usually affect the time evolution of the phenomenon, which are neglected in most PDE control designs. For a vibrating string, the following damping mechanisms are usually neglected:

1) friction with the surrounding medium of the string (or viscous damping; see [12]),

2) thermoelastic phenomena (or thermal damping; see [9,10]),

3) internal friction of the string (or Kelvin-Voigt damping; see [5,17,22,31]).

There is also an additional damping mechanism that may be present when the string has a loose end: a friction mechanism at the loosely attached or damped end. This particular damping mechanism is the so-called "passive damper" and has also been proposed as a boundary feedback controller for the string (see for instance [2,7,8,19,21,24,25]). The passive damper is the simplest boundary feedback law that guarantees exponential stability and it can also guarantee finite-time stability. However, the passive damper has been criticized by many researchers in the literature as a feedback controller that presents severe sensitivity to input delays (see [7,8] as well as the relevant recent discussion in [1] for the case of systems of first-order hyperbolic PDEs). On the other hand, for a "real" string the passive damper is an additional damping mechanism that is (almost) always present.

The present work is devoted to the study of the robustness properties of the 1-D wave equation under all four of the above damping mechanisms. We study robustness with respect to distributed inputs and boundary disturbances in the context of Input-to-State Stability (ISS). The ISS property is a stability property that was proposed by E. D. Sontag in [28] for finite-dimensional control systems and has been recently extended to the case of infinite-dimensional systems described by PDEs (see [16,23] and the references therein). The ISS property and the asymptotic gain property for the wave equation under Kelvin-Voigt damping was recently studied in [17,22].

We consider the 1-D wave equation for a vibrating string with one end pinned and the other end free but in the presence of friction. The same model is obtained if one assumes that the applied force on the free end is manipulated with a passive damper (i.e., using the simplest boundary feedback law) with no input delays. In order to study the robustness properties of the string, we consider different models: (i) a model where only friction with the surrounding medium of the string is taken into account, (ii) a model where thermoelastic phenomena are taken into account as well, and (iii) a model where all damping mechanisms are present. By constructing appropriate ISS Lyapunov functionals for each model, we are in a position to prove the ISS property expressed in different spatial norms.

The structure of the paper is as follows. In Section 2 we present all different mathematical models for a vibrating string. Section 3 of the present work contains the statements of all main results of the paper as well as a discussion on the obtained results. The proofs of the main results are provided in Section 4. Finally, Section 5 gives the concluding remarks of the paper.

**Notation.** Throughout this paper, we adopt the following notation.
* $\Re_+ = [0, +\infty)$ denotes the set of non-negative real numbers.
* Let $S \subseteq \Re^n$ be an open set and let $A \subseteq \Re^n$ be a set that satisfies $S \subseteq A \subseteq cl(S)$. By $C^0(A; \Omega)$, we denote the class of continuous functions on $A$, which take values in $\Omega \subseteq \Re^m$. By $C^k(A; \Omega)$, where $k \geq 1$ is an integer, we denote the class of functions on $A \subseteq \Re^n$, which takes values in $\Omega \subseteq \Re^m$ and has continuous derivatives of order $k$. In other words, the functions of class $C^k(A; \Omega)$ are the functions which have continuous derivatives of order $k$ in $S = \text{int}(A)$ that can be continued continuously to all points in $\partial S \cap A$. When $\Omega = \Re$ then we write $C^0(A)$ or $C^k(A)$. When $I \subseteq \Re$ is an interval and $\eta \in C^1(I)$ is a function of a single variable, $\eta'(\rho)$ denotes the derivative with respect to $\rho \in I$.
* Let $I \subseteq \Re$ be an interval and let $Y$ be a normed linear space. By $C^0(I; Y)$, we denote the class of continuous functions on $I$, which take values in $Y$. By $C^1(I; Y)$, we denote the class of continuously differentiable functions on $I$, which take values in $Y$.
* Let $I \subseteq \Re$ be an interval, let $a < b$ be constants and let $u : I \times [a,b] \to \Re$ be given. We use the notation $u[t]$ to denote the profile at certain $t \in I$, i.e., $(u[t])(x) = u(t, x)$ for all $x \in [a,b]$. When $u(t, x)$ is (twice) differentiable with respect to $x \in [a,b]$, we use the notation $u_x(t,x)$ ($u_{xx}(t,x)$)

for the (second) derivative of $u$ with respect to $x \in [a,b]$, i.e., $u_x(t,x) = \frac{\partial u}{\partial x}(t,x)$ ($u_{xx}(t,x) = \frac{\partial^2 u}{\partial x^2}(t,x)$). When $u(t,x)$ is differentiable with respect to $t$, we use the notation $u_t(t,x)$ for the derivative of $u$ with respect to $t$, i.e., $u_t(t,x) = \frac{\partial u}{\partial t}(t,x)$.

* Let $a < b$ be given constants. For $p \in [1, +\infty)$, $L^p(a,b)$ is the set of equivalence classes of Lebesgue measurable functions $u:(a,b) \to \Re$ with $\|u\|_p := \left( \int_a^b |u(x)|^p \, dx \right)^{1/p} < +\infty$. $L^\infty(a,b)$ is the set of equivalence classes of Lebesgue measurable functions $u:(a,b) \to \Re$ with $\|u\|_\infty := \operatorname*{ess\,sup}_{x \in (a,b)} (|u(x)|) < +\infty$. For an integer $k \geq 1$, $H^k(a,b)$ denotes the Sobolev space of functions in $L^2(a,b)$ with all its weak derivatives up to order $k \geq 1$ in $L^2(a,b)$.

## 2. Four Different Models

The 1-D wave equation

$$u_{tt}(t,x) = c^2 u_{xx}(t,x), \ t > 0, \ x \in (0,1), \tag{2.1}$$

where $c > 0$ is a constant, is a model that is usually employed for the description of the displacement in the $y$-direction $u(t,x)$ at time $t > 0$ and position $x \in (0,1)$ of a string. The model is accompanied by the boundary conditions

$$u(t,0) = 0, \text{ for } t > 0, \tag{2.2}$$

$$u_x(t,1) = U(t), \text{ for } t > 0, \tag{2.3}$$

where $U(t)$ is the control input that corresponds to an external force acting at the right end of the string. The boundary condition (2.2) means that the left end of the string is pinned down.

A family of boundary feedback laws has been proposed in the literature for the global exponential stabilization of the string (see [19]). The feedback laws are given by the formula

$$U(t) = -a u_t(t,1) \text{ -- the passive damper} \tag{2.4}$$

where $a > 0$ is a constant (the controller gain). The feedback law (2.4) (the passive damper) guarantees global exponential stabilization of the string state $(u, u_t)$ in the norm of $H^1(0,1) \times L^2(0,1)$. Moreover, when $a = c^{-1}$ the feedback law (2.4) guarantees finite-time stabilization of the string model (2.1), (2.2), (2.3) (see [2]). It should be noticed that the boundary condition (2.3), (2.4) also arises when the free end of the string moves under friction ("skin friction", modeling the interaction between the surface of the string and the surrounding fluid).

However, as explained in the Introduction, in practice every string is subject to various damping mechanisms and the evolution of the displacement of the string is not described by the simple equation (2.1). When friction with the air (or viscous damping) is taken into account then the closed-loop system is given by (2.2) and

$$u_{tt}(t,x) = c^2 u_{xx}(t,x) - \mu u_t(t,x) + f(t,x), \ t > 0, \ x \in (0,1), \tag{2.5}$$

$$u_x(t,1) = -a u_t(t,1) + d(t), \text{ for } t > 0, \tag{2.6}$$

where $\mu \geq 0$ is a constant (the air friction coefficient-usually unknown), $f$ is a distributed disturbance that may be present and $d$ (a boundary disturbance) can be seen either as the actuator error of the feedback control mechanism or as an unknown external force acting at the end of the string. Model (2.2), (2.5), (2.6) is a model that takes into account viscous damping as well as distributed and boundary disturbances.

A more detailed model is derived when thermoelastic phenomena in the vibrating string are taken into account. In this case the model of the string has an additional state $\theta$ that corresponds to the deviation of the temperature from the reference temperature. Employing the theory of linear thermoelasticity (see [9,10]), we obtain a closed-loop system given by (2.2), (2.6) and the PDEs

$$u_{tt}(t,x) = c^2 u_{xx}(t,x) - \mu u_t(t,x) - b\theta_x(t,x) + f(t,x), \ t > 0, \ x \in (0,1), \tag{2.7}$$

$$\theta_t(t,x) = k\theta_{xx}(t,x) - \lambda u_{xt}(t,x), \ t > 0, \ x \in (0,1), \tag{2.8}$$

where, again, $\mu \geq 0$ is the air friction coefficient, $f$ is a distributed disturbance that may be present, $b, k, \lambda > 0$ are constants that are related to the reference temperature, the reference density of the string, the elastic moduli of the string, the thermal conductivity, the length, and the coefficient of thermal expansion of the string. Assuming that the ends of the string are kept at the constant reference temperature, we also obtain the temperature (Dirichlet) boundary conditions:

$$\theta(t,0) = \theta(t,1) = 0, \text{ for } t > 0. \tag{2.9}$$

Model (2.2), (2.6), (2.7), (2.8), (2.9) is a model that takes into account both viscous and thermal damping as well as distributed and boundary disturbances.

A more complicated model can be obtained if Kelvin-Voigt (internal) damping is taken into account. In this case we replace the PDE (2.7) by the PDE

$$u_{tt}(t,x) = c^2 u_{xx}(t,x) + \sigma u_{xxt}(t,x) - \mu u_t(t,x) - b\theta_x(t,x) + f(t,x), \ t > 0, \ x \in (0,1), \tag{2.10}$$

where $\sigma > 0$ is the coefficient of Kelvin-Voigt damping (a constant). However, in this case, the closed-loop system is not only affected by the actuator error (or external force) $d$ but is also affected by the derivative $\dot{d}$. In order to simplify things, we ignore in this case the control actuator errors (or external force) and consider the model (2.8), (2.10) with boundary conditions (2.2), (2.9) and

$$u_x(t,1) = -a u_t(t,1), \text{ for } t > 0. \tag{2.11}$$

Model (2.2), (2.8), (2.9), (2.10), (2.11) is a model that takes into account viscous, Kelvin-Voigt and thermal damping as well as distributed disturbances. It should be noticed that the PDEs (2.8), (2.10) arise also in the study of a different phenomenon: sound propagation. Indeed, the so-called acoustic approximation for a compressible fluid (linearization around the steady state of constant density and temperature and zero fluid velocity) moving in a single spatial direction gives the PDEs (equations of viscous thermoacoustics-see [18]):

$$\rho_t + \gamma v_x = 0 \tag{2.12}$$

$$v_t + \frac{c^2}{\gamma}\rho_x + b\theta_x = \sigma v_{xx} \tag{2.13}$$

$$\theta_t + \lambda v_x = k\theta_{xx} \tag{2.14}$$

where $c, \gamma, b, k, \lambda, \sigma > 0$ are constants (depending on the physical properties of the fluid at the reference density and reference temperature), $v$ is the fluid velocity and $\rho, \theta$ are the deviations of density, temperature, respectively, from their reference values. Differentiating (2.12) with respect to $x$ and (2.13) with respect to $t$ we get the PDE:

$$v_{tt} = c^2 v_{xx} + \sigma v_{xxt} - b\theta_{xt} \tag{2.15}$$

Setting $v = u_t$, the PDEs (2.14), (2.15) give us the PDEs (2.8), (2.10) with $\mu = 0$ (since there is no surrounding medium like air, there is no friction with the surrounding medium).

Table 1 shows the four different models for a vibrating string under the feedback law (2.4).

**Table 1:** Four models for a vibrating string under the feedback law (2.4)

| Closed-Loop Models | Equations | Viscous damping | Thermal damping | Kelvin-Voigt damping | Distributed disturbance | Control Actuator Error |
|---|---|---|---|---|---|---|
| Model (A) | (2.1), (2.2), (2.11) | - | - | - | - | - |
| Model (B) | (2.2), (2.5), (2.6) | ✓ | - | - | ✓ | ✓ |
| Model (C) | (2.2), (2.6), (2.7), (2.8), (2.9) | ✓ | ✓ | - | ✓ | ✓ |
| Model (D) | (2.2), (2.8), (2.9), (2.10), (2.11) | ✓ | ✓ | ✓ | ✓ | - |

The feedback law (2.4) is designed for the simplest string model (2.1), (2.2), (2.3). Moreover, it is generally difficult to know the exact values of all constants $b, k, \lambda > 0$, $\sigma, \mu \geq 0$ that appear in the more complicated models. A boundary feedback design procedure that is based on the more complicated models would require knowledge of the exact values of all constants $b, k, \lambda > 0$, $\sigma, \mu \geq 0$. Consequently, two important robustness questions arise for the family of simple feedback laws (2.4):

1) Does the string exhibit exponential stability in the absence of disturbances and in some appropriate state norm for Model (B), Model (C) and Model (D)?

2) Does the string satisfy an ISS property with respect to the disturbances and in some appropriate state norm for Model (B), Model (C) and Model (D)?

The next section of the paper is devoted to the answers of these questions.

## 3. Main Results

The answers to the questions that were posed in the previous section are positive. In other words, for all models the string exhibits exponential stability in the absence of disturbances and for all models the string satisfies an ISS property with respect to the disturbances. The following theorems give precise information for each model and are the main results of the present work.

**Theorem 1 (String with Viscous Damping):** *For every $a,c > 0$, $\mu \geq 0$ there exist constants $\omega, G, \gamma_1, \gamma_2 > 0$ such that every solution $u \in C^1(\Re_+ \times [0,1]) \cap C^2((0,+\infty) \times (0,1))$ of (2.2), (2.5), (2.6) corresponding to inputs $f \in C^0(\Re_+; L^2(0,1)) \cap C^0((0,+\infty) \times (0,1))$, $d \in C^0(\Re_+)$ satisfies the following estimate for all $t \geq 0$:*

$$\left(\|u_t[t]\|_2^2 + \|u_x[t]\|_2^2\right)^{1/2} \leq G\exp(-\omega t)\left(\|u_t[0]\|_2^2 + \|u_x[0]\|_2^2\right)^{1/2} + \gamma_1 \sup_{0 \leq s \leq t}\left(\|f[s]\|_2\right) + \gamma_2 \sup_{0 \leq s \leq t}\left(|d(s)|\right) \quad (3.1)$$

**Theorem 2 (String with Viscous and Thermal Damping):** *For every $a,c > 0$, $\mu \geq 0$, $b,k,\lambda > 0$ there exist constants $\omega, G, \gamma_1, \gamma_2 > 0$ such that every solution $u \in C^1(\Re_+ \times [0,1]) \cap C^2((0,+\infty) \times (0,1))$, $\theta \in C^0(\Re_+ \times [0,1]) \cap C^1((0,+\infty); L^2(0,1))$ with $\theta[t] \in H^2(0,1)$ for all $t > 0$ of (2.2), (2.6), (2.7), (2.8), (2.9) corresponding to inputs $f \in C^0(\Re_+; L^2(0,1)) \cap C^0((0,+\infty) \times (0,1))$, $d \in C^0(\Re_+)$ satisfies the following estimate for all $t \geq 0$:*

$$\left(\|u_t[t]\|_2^2 + \|u_x[t]\|_2^2 + \|\theta[t]\|_2^2\right)^{1/2} \leq G\exp(-\omega t)\left(\|u_t[0]\|_2^2 + \|u_x[0]\|_2^2 + \|\theta[0]\|_2^2\right)^{1/2} + \gamma_1 \sup_{0 \leq s \leq t}\left(\|f[s]\|_2\right) + \gamma_2 \sup_{0 \leq s \leq t}\left(|d(s)|\right) \quad (3.2)$$

**Theorem 3 (String with Viscous, Thermal, and Kelvin-Voigt Damping):** *For every $a,c > 0$, $\mu \geq 0$, $\sigma, b, k, \lambda > 0$ there exist constants $\omega, G, \gamma > 0$ such that every solution $u \in C^2(\Re_+ \times [0,1])$, $\theta \in C^0(\Re_+ \times [0,1]) \cap C^1((0,+\infty); L^2(0,1))$ with $u_{xx} \in C^1((0,+\infty); L^2(0,1))$, $\theta[t], u_t[t] \in H^2(0,1)$ for all $t > 0$ of (2.2), (2.8), (2.9), (2.10), (2.11) corresponding to input $f \in C^0(\Re_+; L^2(0,1)) \cap C^0((0,+\infty) \times (0,1))$ satisfies the following estimate for all $t \geq 0$:*

$$\left(\|u_t[t]\|_2^2 + \|u_x[t]\|_2^2 + \|u_{xx}[t]\|_2^2 + |u_t(t,1)|^2 + \|\theta[t]\|_2^2\right)^{1/2}$$
$$\leq G\exp(-\omega t)\left(\|u_t[0]\|_2^2 + \|u_x[0]\|_2^2 + \|u_{xx}[0]\|_2^2 + |u_t(0,1)|^2 + \|\theta[0]\|_2^2\right)^{1/2} + \gamma \sup_{0 \leq s \leq t}\left(\|f[s]\|_2\right) \quad (3.3)$$

**Remarks: (i)** The reader should notice that ISS estimates (3.1), (3.2), (3.3) for each of the models (model (B), model (C) and model (D); recall Table 1) actually guarantee the so-called "exponential ISS" property (see [23]), i.e., the effect of the initial condition is expressed by an exponential *KL* function and the gain functions of the disturbances are linear. The ISS estimates (3.1), (3.2), (3.3) show global exponential stability in the absence of disturbances.

**(ii)** The ISS estimates (3.1), (3.2), (3.3) for each of the models (model (B), model (C) and model (D); recall Table 1) are expressed in different state spatial norms. More specifically, ISS estimate (3.1) for model (B) involves the norm of the space $H^1(0,1) \times L^2(0,1)$, ISS estimate (3.2) for model (C) involves the norm of the space $H^1(0,1) \times L^2(0,1) \times L^2(0,1)$, while ISS estimate (3.3) for model (D) involves a norm of the space $H^2(0,1) \times C^0([0,1]) \times L^2(0,1)$ where the norm of the state component $u_t \in C^0([0,1])$ is not the standard sup-norm but the norm $\|u_t\| = \left( \|u_t\|_2^2 + |u_t(1)|^2 \right)^{1/2}$.

**(iii)** Each of the above theorems (Theorem 1, Theorem 2 and Theorem 3) assumes different regularity properties for the solution (going from less regular solutions to more regular solutions). This is expected, since the models (model (B), model (C) and model (D); recall Table 1) are quite different models involving very different differential operators. The regularity requirements for the solutions can be guaranteed by assuming sufficiently regular initial conditions and sufficiently regular applied inputs that also satisfy additional compatibility requirements. However, the specification of the (sufficient) regularity requirements for the initial conditions and the applied inputs is out of the scope of the present paper.

**(iv)** Each of the above theorems (Theorem 1, Theorem 2 and Theorem 3) provide qualitative results for the solutions of the corresponding closed-loop systems. In other words, the above theorems do not provide formulae that show how large are the constants that are involved in the corresponding ISS estimates. However, it should be noticed that the proofs of the above theorems are constructive and provide estimates for the constants.

**(v)** The proofs of the above theorems follow a similar methodology: the construction of an ISS Lyapunov functional. However, for each of the closed-loop systems (model (B), model (C) and model (D); recall Table 1) the ISS Lyapunov functionals are different.

**(vii)** As pointed out above, when $a = c^{-1}$ the feedback law (2.4) guarantees finite-time stabilization of the string model (2.1), (2.2), (2.3) (see [2]). However, this is not the case when the closed-loop system is described by one of the more complicated models (model (B), model (C) and model (D); recall Table 1): finite-time stability is a feature that is not preserved under perturbations of the original string model.

## 4. Proofs of Main Results

We start by providing the proof of Theorem 1.

**Proof of Theorem 1:** Consider an arbitrary solution $u \in C^1(\Re_+ \times [0,1]) \cap C^2((0,+\infty) \times (0,1))$ of (2.2), (2.5), (2.6) corresponding to (arbitrary) inputs $f \in C^0(\Re_+; L^2(0,1)) \cap C^0((0,+\infty) \times (0,1))$, $d \in C^0(\Re_+)$. Let $r > 0$ be a constant and define the functionals $E, \Phi : H^1(0,1) \times L^2(0,1) \to \Re_+$ by means of the following formulae:

$$E(u,w) := \frac{1}{2}\int_0^1 w^2(x)dx + \frac{c^2}{2}\int_0^1 u_x^2(x)dx, \text{ for all } u \in H^1(0,1), w \in L^2(0,1) \quad (4.1)$$

$$\Phi(u,w) := \frac{1}{2}\int_0^1 \exp(rx)(w(x) + cu_x(x))^2 dx + \frac{1}{2}\int_0^1 \exp(-rx)(w(x) - cu_x(x))^2 dx,$$
$$\text{for all } u \in H^1(0,1), w \in L^2(0,1) \quad (4.2)$$

It should be noticed at this point that $E(u[t], u_t[t])$ is the mechanical energy of the string at time $t \geq 0$, while $\Phi(u[t], u_t[t])$ is the value of a Lyapunov functional for the simple model (A). When disturbances are absent and when $\mu = 0$ (no viscous damping), the function $\Phi(u[t], u_t[t])$ can be used for the derivation of an exponential stability estimate for the string.

Since $r > 0$, the functional $\Phi : H^1(0,1) \times L^2(0,1) \to \mathfrak{R}_+$ satisfies the following estimate:

$$\exp(-r)\left(c^2 \|u_x\|_2^2 + \|w\|_2^2\right) \leq \Phi(u,w) \leq \exp(r)\left(c^2 \|u_x\|_2^2 + \|w\|_2^2\right), \text{ for all } u \in H^1(0,1),\ w \in L^2(0,1) \quad (4.3)$$

Since $u \in C^2\left((0,+\infty) \times (0,1)\right)$, we get from definitions (4.1), (4.2) for all $t > 0$:

$$\frac{d}{dt} E(u[t], u_t[t]) = \int_0^1 u_t(t,x) u_{tt}(t,x) dx + c^2 \int_0^1 u_x(t,x) u_{xt}(t,x) dx \quad (4.4)$$

$$\frac{d}{dt} \Phi(u[t], u_t[t]) = \int_0^1 \exp(rx)\left(u_t(t,x) + c u_x(t,x)\right)\left(u_{tt}(t,x) + c u_{xt}(t,x)\right) dx$$
$$+ \int_0^1 \exp(-rx)\left(u_t(t,x) - c u_x(t,x)\right)\left(u_{tt}(t,x) - c u_{xt}(t,x)\right) dx \quad (4.5)$$

Integrating by parts the integral $\int_0^1 u_x(t,x) u_{xt}(t,x) dx$ and using (4.4), (2.5), (2.6), (2.2) (the latter implies that $u_t(t,0) = 0$ for all $t > 0$) we obtain for $t > 0$:

$$\frac{d}{dt} E(u[t], u_t[t]) = -\mu \int_0^1 u_t^2(t,x) dx + \int_0^1 u_t(t,x) f(t,x) dx - ac^2 u_t^2(t,1) + c^2 d(t) u_t(t,1) \quad (4.6)$$

Using (4.5) and (2.5) we get for $t > 0$:

$$\frac{d}{dt} \Phi(u[t], u_t[t]) = c \int_0^1 \exp(rx) \left(\frac{1}{2}\left(u_t(t,x) + c u_x(t,x)\right)^2\right)_x dx$$
$$- c \int_0^1 \exp(-rx) \left(\frac{1}{2}\left(u_t(t,x) - c u_x(t,x)\right)^2\right)_x dx$$
$$+ \int_0^1 \exp(rx)\left(u_t(t,x) + c u_x(t,x)\right) f(t,x) dx$$
$$+ \int_0^1 \exp(-rx)\left(u_t(t,x) - c u_x(t,x)\right) f(t,x) dx \quad (4.7)$$
$$- \mu \int_0^1 \exp(rx)\left(u_t(t,x) + c u_x(t,x)\right) u_t(t,x) dx$$
$$- \mu \int_0^1 \exp(-rx)\left(u_t(t,x) - c u_x(t,x)\right) u_t(t,x) dx$$

Using the inequalities $d(t) u_t(t,1) \leq \frac{a}{2} u_t^2(t,1) + \frac{1}{2a} |d(t)|^2$ (which holds since $a > 0$),

$$\left(u_t(t,x) \pm c u_x(t,x)\right) f(t,x) \leq \frac{cr}{4(1+\mu)} \left(u_t(t,x) \pm c u_x(t,x)\right)^2 + \frac{1+\mu}{cr} |f(t,x)|^2,$$

$$\left(u_t(t,x) \pm c u_x(t,x)\right) u_t(t,x) \geq -\frac{cr}{4(1+\mu)} \left(u_t(t,x) \pm c u_x(t,x)\right)^2 - \frac{1+\mu}{cr} u_t^2(t,x)$$

and the fact that $\mu \geq 0$, we obtain from (4.6), (4.7) for $t > 0$:

$$\frac{d}{dt} E(u[t], u_t[t]) \leq -\mu \int_0^1 u_t^2(t,x) dx + \int_0^1 u_t(t,x) f(t,x) dx - \frac{ac^2}{2} u_t^2(t,1) + \frac{c^2}{2a} |d(t)|^2 \qquad (4.8)$$

$$\frac{d}{dt} \Phi(u[t], u_t[t]) \leq c \int_0^1 \exp(rx) \left(\frac{1}{2} \left(u_t(t,x) + c u_x(t,x)\right)^2\right)_x dx$$

$$-c \int_0^1 \exp(-rx) \left(\frac{1}{2} \left(u_t(t,x) - c u_x(t,x)\right)^2\right)_x dx$$

$$+ \frac{cr}{4} \int_0^1 \exp(rx) \left(u_t(t,x) + c u_x(t,x)\right)^2 dx \qquad (4.9)$$

$$+ \frac{cr}{4} \int_0^1 \exp(-rx) \left(u_t(t,x) - c u_x(t,x)\right)^2 dx$$

$$+ \frac{2\mu(1+\mu)}{cr} \int_0^1 \cosh(rx) u_t^2(t,x) dx + \frac{2(1+\mu)}{cr} \int_0^1 \cosh(rx) |f(t,x)|^2 dx$$

Integrating by parts the integrals $\int_0^1 \exp(rx) \left(\frac{1}{2} \left(u_t(t,x) + c u_x(t,x)\right)^2\right)_x dx$, $\int_0^1 \exp(-rx) \left(\frac{1}{2} \left(u_t(t,x) - c u_x(t,x)\right)^2\right)_x dx$ and using (4.9), (2.6), (2.2) (the latter implies that $u_t(t,0) = 0$ for all $t > 0$) we obtain for $t > 0$:

$$\frac{d}{dt} \Phi(u[t], u_t[t])$$

$$\leq \frac{c}{2} \left(\exp(r) \left((1-ac) u_t(t,1) + cd(t)\right)^2 - \exp(-r) \left((1+ac) u_t(t,1) - cd(t)\right)^2\right)$$

$$- \frac{cr}{4} \int_0^1 \exp(rx) \left(u_t(t,x) + c u_x(t,x)\right)^2 dx - \frac{cr}{4} \int_0^1 \exp(-rx) \left(u_t(t,x) - c u_x(t,x)\right)^2 dx \qquad (4.10)$$

$$+ \frac{2\mu(1+\mu)}{cr} \int_0^1 \cosh(rx) u_t^2(t,x) dx + \frac{2(1+\mu)}{cr} \int_0^1 \cosh(rx) |f(t,x)|^2 dx$$

The fact that $2 u_t(t,x) = \left(u_t(t,x) + c u_x(t,x)\right) + \left(u_t(t,x) - c u_x(t,x)\right)$ implies that

$$4 u_t^2(t,x) \leq (1+\zeta) \left(u_t(t,x) + c u_x(t,x)\right)^2 + (1+\zeta^{-1}) \left(u_t(t,x) - c u_x(t,x)\right)^2, \text{ for all } \zeta > 0$$

Setting $\zeta = \exp(2rx)$ we get

$$2 u_t^2(t,x) \leq \cosh(rx) \left(\exp(rx) \left(u_t(t,x) + c u_x(t,x)\right)^2 + \exp(-rx) \left(u_t(t,x) - c u_x(t,x)\right)^2\right)$$

Consequently, we get:

$$\int_0^1 u_t^2(t,x)dx \leq \frac{\cosh(r)}{2}\int_0^1 \exp(rx)\left(u_t(t,x)+cu_x(t,x)\right)^2 dx$$
$$+\frac{\cosh(r)}{2}\int_0^1 \exp(-rx)\left(u_t(t,x)-cu_x(t,x)\right)^2 dx \tag{4.11}$$

Define

$$M := \max\left(\frac{2(1+\mu)}{cr}\cosh(r), \frac{2}{ac}\left(\exp(r)(1-ac)^2 - \exp(-r)(1+ac)^2\right)\right) \tag{4.12}$$

Using the inequality $u_t(t,x)f(t,x) \leq \frac{cr}{4M\cosh(r)}u_t^2(t,x) + \frac{M\cosh(r)}{cr}|f(t,x)|^2$ and (4.8), (4.11) we obtain for $t>0$:

$$\frac{d}{dt}E(u[t],u_t[t]) \leq -\mu\int_0^1 u_t^2(t,x)dx + \frac{cr}{8M}\int_0^1 \exp(rx)\left(u_t(t,x)+cu_x(t,x)\right)^2 dx$$
$$+\frac{cr}{8M}\int_0^1 \exp(-rx)\left(u_t(t,x)-cu_x(t,x)\right)^2 dx \tag{4.13}$$
$$-\frac{ac^2}{2}u_t^2(t,1) + \frac{c^2}{2a}|d(t)|^2 + \frac{M\cosh(r)}{cr}\|f[t]\|_2^2$$

Define the functional $V: H^1(0,1)\times L^2(0,1) \to \Re_+$ by means of the formula

$$V(u,w) := \Phi(u,w) + ME(u,w), \tag{4.14}$$

with $M$ given in (4.12). Using (4.1), (4.2), (4.3), (4.10), (4.13) and definition (4.14) we conclude that the functional $V: H^1(0,1)\times L^2(0,1) \to \Re_+$ satisfies the following estimates:

$$\left(\frac{M}{2}+\exp(-r)\right)\left(c^2\|u_x\|_2^2+\|w\|_2^2\right) \leq V(u,w) \leq \left(\frac{M}{2}+\exp(r)\right)\left(c^2\|u_x\|_2^2+\|w\|_2^2\right),$$
$$\text{for all } u \in H^1(0,1),\ w \in L^2(0,1) \tag{4.15}$$

$$\frac{d}{dt}V(u[t],u_t[t]) \leq -\frac{c}{2}\left(acM + \exp(-r)(1+ac)^2 - \exp(r)(1-ac)^2\right)u_t^2(t,1)$$
$$+2c^2\left(\cosh(r)-ac\sinh(r)\right)u_t(t,1)d(t) + c^2\left(\frac{M}{2a}+\sinh(r)c\right)|d(t)|^2 - \frac{cr}{4}\Phi(u[t],u_t[t])$$
$$-\mu\left(M - \frac{2(1+\mu)}{cr}\cosh(r)\right)\int_0^1 u_t^2(t,x)dx + \left(2(1+\mu)+M^2\right)\frac{\cosh(r)}{cr}\|f[t]\|_2^2$$
$$\text{for all } t>0 \tag{4.16}$$

Using the facts that $\mu \geq 0$, $M \geq \frac{2(1+\mu)}{cr}\cosh(r)$, $M \geq \frac{2}{ac}\left(\exp(r)(1-ac)^2 - \exp(-r)(1+ac)^2\right)$ (recall (4.12)) and the inequality

$$\left(\cosh(r)-ac\sinh(r)\right)u_t(t,1)d(t) \leq \frac{aM}{8}u_t^2(t,1) + \frac{2\left(\cosh(r)-ac\sinh(r)\right)^2}{aM}|d(t)|^2$$

we conclude from (4.16) that the following estimate holds for all $t>0$:

$$\frac{d}{dt}V(u[t],u_t[t]) \leq -\frac{cr}{4}\Phi\left(u[t],u_t[t]\right) + \left(2(1+\mu)+M^2\right)\frac{\cosh(r)}{cr}\|f[t]\|_2^2$$
$$+c^2\left(\frac{M}{2a} + \frac{4\left(\cosh(r)-ac\sinh(r)\right)^2}{aM} + c\sinh(r)\right)|d(t)|^2 \qquad (4.17)$$

Definitions (4.1), (4.14) and estimate (4.15) imply that

$$-\Phi(u,w) = -V(u,w)+ME(u,w) \leq -V(u,w)+\frac{M}{M+2\exp(-r)}V(u,w) = -\frac{2\exp(-r)}{M+2\exp(-r)}V(u,w)$$

for all $u \in H^1(0,1)$, $w \in L^2(0,1)$. The above inequality implies the following differential inequality for all $t>0$:

$$\frac{d}{dt}V(u[t],u_t[t]) \leq -\frac{cr\exp(-r)}{2(M+2\exp(-r))}V\left(u[t],u_t[t]\right) + \left(2(1+\mu)+M^2\right)\frac{\cosh(r)}{cr}\|f[t]\|_2^2$$
$$+c^2\left(\frac{M}{2a} + \frac{4\left(\cosh(r)-ac\sinh(r)\right)^2}{aM} + c\sinh(r)\right)|d(t)|^2 \qquad (4.18)$$

Differential inequality (4.18) directly implies the following estimate for all $t \geq t_0 > 0$:

$$V(u[t],u_t[t]) \leq \exp\left(-2\omega(t-t_0)\right)V\left(u[t_0],u_t[t_0]\right) + K_1 \sup_{t_0 \leq s \leq t}\left(\|f[s]\|_2^2\right) + K_2 \sup_{t_0 \leq s \leq t}\left(|d(s)|^2\right) \qquad (4.19)$$

with $\omega := \dfrac{cr\exp(-r)}{4(M+2\exp(-r))}$, $K_1 := \left(2(1+\mu)+M^2\right)\dfrac{\cosh(r)}{2\omega cr}$ and

$K_2 := \dfrac{c^2}{2\omega}\left(\dfrac{M}{2a} + \dfrac{4\left(\cosh(r)-ac\sinh(r)\right)^2}{aM} + c\sinh(r)\right)$. Since $u \in C^1\left(\Re_+ \times [0,1]\right)$ we obtain that $\lim_{t_0 \to 0^+}\left(V\left(u[t_0],u_t[t_0]\right)\right) = V\left(u[0],u_t[0]\right)$ and consequently, we get the following estimate for all $t \geq 0$:

$$V(u[t],u_t[t]) \leq \exp(-2\omega t)V\left(u[0],u_t[0]\right) + K_1 \sup_{0 \leq s \leq t}\left(\|f[s]\|_2^2\right) + K_2 \sup_{0 \leq s \leq t}\left(|d(s)|^2\right) \qquad (4.20)$$

Exploiting estimates (4.20) and (4.15) we get for all $t \geq 0$:

$$c^2\|u_x[t]\|_2^2 + \|u_t[t]\|_2^2 \leq \exp(-2\omega t)\frac{M+2\exp(r)}{M+2\exp(-r)}\left(c^2\|u_x[0]\|_2^2 + \|u_t[0]\|_2^2\right)$$
$$+\frac{2K_1}{M+2\exp(-r)}\sup_{0 \leq s \leq t}\left(\|f[s]\|_2^2\right) + \frac{2K_2}{M+2\exp(-r)}\sup_{0 \leq s \leq t}\left(|d(s)|^2\right) \qquad (4.21)$$

Estimate (3.1) with $G := \sqrt{\dfrac{M+2\exp(r)}{M+2\exp(-r)}} \sqrt{\dfrac{\max(1,c^2)}{\min(1,c^2)}}$, $\gamma_1 := \sqrt{\dfrac{2K_1}{(M+2\exp(-r))\min(1,c^2)}}$ and $\gamma_2 := \gamma_1 \sqrt{\dfrac{K_2}{K_1}}$ is a direct consequence of estimate (4.21). The proof is complete. ◁

We next continue with the proof of Theorem 2.

**Proof of Theorem 2:** The proof of Theorem 2 uses a similar notation with the proof of Theorem 1 (for example, there is a constant $M > 0$, there are functionals $E, V$, etc.). However, the reader should not be tempted, by an overlapping notation for different quantities, to compare the different quantities in the proofs of Theorem 1 and Theorem 2. To understand this, we point out that while the functional $E$ in the proof of Theorem 1 is the mechanical energy of the string, the functional $E$ in the proof of Theorem 2 below is not the mechanical energy of the string and is not a functional that may be considered as some kind of energy of the string (e.g., the sum of mechanical energy and the thermal energy of the string).

Consider an arbitrary solution $u \in C^1(\Re_+ \times [0,1]) \cap C^2((0,+\infty) \times (0,1))$, $\theta \in C^0(\Re_+ \times [0,1]) \cap C^1((0,+\infty); L^2(0,1))$ with $\theta[t] \in H^2(0,1)$ for all $t > 0$ of (2.2), (2.6), (2.7), (2.8), (2.9) corresponding to (arbitrary) inputs $f \in C^0(\Re_+; L^2(0,1)) \cap C^0((0,+\infty) \times (0,1))$, $d \in C^0(\Re_+)$. Let $r > 0$ be a constant and define the functionals $E : H^1(0,1) \times L^2(0,1) \times L^2(0,1) \to \Re_+$, $\Phi : H^1(0,1) \times L^2(0,1) \to \Re_+$ by means of the formulae (4.2) and

$$E(u,w,\theta) := \frac{1}{2}\int_0^1 w^2(x)dx + \frac{c^2}{2}\int_0^1 u_x^2(x)dx + \frac{b}{2\lambda}\int_0^1 \theta^2(x)dx, \text{ for all } u \in H^1(0,1),\ w,\theta \in L^2(0,1)$$

(4.22)

It should be noticed at this point that $E(u[t], u_t[t], \theta[t])$ is an appropriate linear combination at time $t \geq 0$ of the mechanical energy ($\frac{1}{2}\int_0^1 u_t^2(t,x)dx + \frac{c^2}{2}\int_0^1 u_x^2(t,x)dx$) of the string and the squared $L^2$ norm of the deviation of temperature from its reference value (the term $\int_0^1 \theta^2(t,x)dx$ is not the thermal energy of the string but is the value of a Lyapunov functional used frequently for the heat equation).

Since $r > 0$, the functional $\Phi : H^1(0,1) \times L^2(0,1) \to \Re_+$ satisfies estimate (4.3). Moreover, since $u \in C^2((0,+\infty) \times (0,1))$, $\theta \in C^1((0,+\infty); L^2(0,1))$ we get from definitions (4.2), (4.22) for all $t > 0$ formula (4.5) as well as the following formula:

$$\frac{d}{dt}E(u[t], u_t[t], \theta[t]) = \int_0^1 u_t(t,x)u_{tt}(t,x)dx + c^2\int_0^1 u_x(t,x)u_{xt}(t,x)dx + \frac{b}{\lambda}\int_0^1 \theta(t,x)\theta_t(t,x)dx \quad (4.23)$$

Integrating by parts the integral $\int_0^1 u_x(t,x)u_{xt}(t,x)dx$ and using (4.23), (2.6), (2.7), (2.8), (2.2) (the latter implies that $u_t(t,0) = 0$ for all $t > 0$) we obtain for $t > 0$:

$$\frac{d}{dt}E(u[t],u_t[t],\theta[t]) = -\mu\int_0^1 u_t^2(t,x)dx + \int_0^1 u_t(t,x)f(t,x)dx - ac^2 u_t^2(t,1)$$
$$+c^2 d(t)u_t(t,1) - b\int_0^1 \left(u_t(t,x)\theta(t,x)\right)_x dx + \frac{kb}{\lambda}\int_0^1 \theta(t,x)\theta_{xx}(t,x)dx \tag{4.24}$$

Using the inequality $d(t)u_t(t,1) \leq \frac{a}{2}u_t^2(t,1) + \frac{1}{2a}|d(t)|^2$ (which holds since $a>0$), (2.9) and integrating by parts the integral $\int_0^1 \theta(t,x)\theta_{xx}(t,x)dx$, we obtain from (4.24) for $t>0$:

$$\frac{d}{dt}E(u[t],u_t[t],\theta[t]) \leq -\mu\int_0^1 u_t^2(t,x)dx + \int_0^1 u_t(t,x)f(t,x)dx$$
$$-\frac{ac^2}{2}u_t^2(t,1) + \frac{c^2}{2a}|d(t)|^2 - \frac{kb}{\lambda}\int_0^1 \theta_x^2(t,x)dx \tag{4.25}$$

Using (4.5) and (2.7) we get for $t>0$:

$$\frac{d}{dt}\Phi(u[t],u_t[t]) = c\int_0^1 \exp(rx)\left(\frac{1}{2}\left(u_t(t,x)+cu_x(t,x)\right)^2\right)_x dx$$
$$-c\int_0^1 \exp(-rx)\left(\frac{1}{2}\left(u_t(t,x)-cu_x(t,x)\right)^2\right)_x dx$$
$$+\int_0^1 \exp(rx)\left(u_t(t,x)+cu_x(t,x)\right)\left(f(t,x)-b\theta_x(t,x)-\mu u_t(t,x)\right)dx \tag{4.26}$$
$$+\int_0^1 \exp(-rx)\left(u_t(t,x)-cu_x(t,x)\right)\left(f(t,x)-b\theta_x(t,x)-\mu u_t(t,x)\right)dx$$

Integrating by parts the integrals $\int_0^1 \exp(rx)\left(\frac{1}{2}\left(u_t(t,x)+cu_x(t,x)\right)^2\right)_x dx$, $\int_0^1 \exp(-rx)\left(\frac{1}{2}\left(u_t(t,x)-cu_x(t,x)\right)^2\right)_x dx$ and using (4.26), (2.6), (2.2) (the latter implies that $u_t(t,0)=0$ for all $t>0$) we obtain for $t>0$:

$$\frac{d}{dt}\Phi(u[t],u_t[t]) = -\frac{cr}{2}\int_0^1 \exp(rx)\left(u_t(t,x)+cu_x(t,x)\right)^2 dx$$
$$-\frac{cr}{2}\int_0^1 \exp(-rx)\left(u_t(t,x)-cu_x(t,x)\right)^2 dx$$
$$+\frac{c}{2}\left(\exp(r)\left((1-ac)u_t(t,1)+cd(t)\right)^2 - \exp(-r)\left((1+ac)u_t(t,1)-cd(t)\right)^2\right) \tag{4.27}$$
$$+\int_0^1 \exp(rx)\left(u_t(t,x)+cu_x(t,x)\right)\left(f(t,x)-b\theta_x(t,x)-\mu u_t(t,x)\right)dx$$
$$+\int_0^1 \exp(-rx)\left(u_t(t,x)-cu_x(t,x)\right)\left(f(t,x)-b\theta_x(t,x)-\mu u_t(t,x)\right)dx$$

Using the inequalities

$$\left(u_t(t,x) \pm c u_x(t,x)\right) f(t,x) \le \frac{cr}{4(1+b+\mu)} \left(u_t(t,x) \pm c u_x(t,x)\right)^2 + \frac{1+b+\mu}{cr} |f(t,x)|^2,$$

$$\left(u_t(t,x) \pm c u_x(t,x)\right) u_t(t,x) \ge -\frac{cr}{4(1+b+\mu)} \left(u_t(t,x) \pm c u_x(t,x)\right)^2 - \frac{1+b+\mu}{cr} u_t^2(t,x)$$

$$\left(u_t(t,x) \pm c u_x(t,x)\right) u_t(t,x) \ge -\frac{cr}{4(1+b+\mu)} \left(u_t(t,x) \pm c u_x(t,x)\right)^2 - \frac{1+b+\mu}{cr} u_t^2(t,x)$$

and the facts that $\mu \ge 0$, $b > 0$, we obtain from (4.27) the following estimate for $t > 0$:

$$\begin{aligned}
\frac{d}{dt} \Phi(u[t], u_t[t]) &\le -\frac{cr}{4} \int_0^1 \exp(rx) \left(u_t(t,x) + c u_x(t,x)\right)^2 dx \\
&\quad -\frac{cr}{4} \int_0^1 \exp(-rx) \left(u_t(t,x) - c u_x(t,x)\right)^2 dx \\
&\quad -\frac{c}{2} \left(\exp(-r)(1+ac)^2 - \exp(r)(1-ac)^2\right) u_t^2(t,1) \\
&\quad + \sinh(r) c^3 |d(t)|^2 + 2c^2 \left(\cosh(r) - ac \sinh(r)\right) u_t(t,1) d(t) \\
&\quad + \frac{2(1+b+\mu)}{cr} \int_0^1 \cosh(rx) |f(t,x)|^2 dx + \frac{2b(1+b+\mu)}{cr} \int_0^1 \cosh(rx) \theta_x^2(t,x) dx \\
&\quad + \frac{2\mu(1+b+\mu)}{cr} \int_0^1 \cosh(rx) u_t^2(t,x) dx
\end{aligned} \quad (4.28)$$

Define

$$M := \max\left(\max\left(1, \frac{2\lambda}{k}\right) \frac{2(1+b+\mu)}{cr} \cosh(r), \frac{2}{ac}\left(\exp(r)(1-ac)^2 - \exp(-r)(1+ac)^2\right)\right) \quad (4.29)$$

$$V(u, w, \theta) := \Phi(u, w) + M E(u, w, \theta), \text{ for all } u \in H^1(0,1), \ w, \theta \in L^2(0,1) \quad (4.30)$$

Using (4.25), (4.28) and definition (4.30), we get for $t > 0$:

$$\begin{aligned}
\frac{d}{dt} V(u[t], u_t[t], \theta[t]) &\le -\frac{cr}{4} \int_0^1 \exp(rx) \left(u_t(t,x) + c u_x(t,x)\right)^2 dx \\
&\quad -\frac{cr}{4} \int_0^1 \exp(-rx) \left(u_t(t,x) - c u_x(t,x)\right)^2 dx + \frac{2(1+b+\mu)}{cr} \cosh(r) \|f[t]\|_2^2 \\
&\quad -\frac{c}{2} \left(acM + \exp(-r)(1+ac)^2 - \exp(r)(1-ac)^2\right) u_t^2(t,1) \\
&\quad + c^2 \left(c \sinh(r) + \frac{M}{2a}\right) |d(t)|^2 + 2c^2 \left(\cosh(r) - ac \sinh(r)\right) u_t(t,1) d(t) \\
&\quad -\mu \left(M - \frac{2(1+b+\mu)}{cr} \cosh(r)\right) \int_0^1 u_t^2(t,x) dx + M \int_0^1 u_t(t,x) f(t,x) dx \\
&\quad -b\left(\frac{kM}{\lambda} - \frac{2(1+b+\mu)}{cr} \cosh(r)\right) \int_0^1 \theta_x^2(t,x) dx
\end{aligned} \quad (4.31)$$

Using the inequality $\int_0^1 u_t(t,x) f(t,x) dx \leq \dfrac{cr}{4M \cosh(r)} \int_0^1 u_t^2(t,x) dx + \dfrac{M \cosh(r)}{cr} \|f[t]\|_2^2$, inequality (4.11), the fact that $M \geq \dfrac{2(1+b+\mu)}{cr} \cosh(r)$ (recall (4.29)), the fact that $\dfrac{kM}{2\lambda} \geq \dfrac{2(1+b+\mu)}{cr} \cosh(r)$ (recall (4.29)) and the facts that $\mu \geq 0$, $b > 0$, we conclude from (4.31) the following estimate for $t > 0$:

$$\begin{aligned}
\frac{d}{dt} V(u[t], u_t[t], \theta[t]) &\leq -\frac{cr}{8} \int_0^1 \exp(rx) \left( u_t(t,x) + cu_x(t,x) \right)^2 dx \\
&\quad -\frac{cr}{8} \int_0^1 \exp(-rx) \left( u_t(t,x) - cu_x(t,x) \right)^2 dx \\
&\quad -\frac{c}{2} \left( acM + \exp(-r)(1+ac)^2 - \exp(r)(1-ac)^2 \right) u_t^2(t,1) \\
&\quad + c^2 \left( c \sinh(r) + \frac{M}{2a} \right) |d(t)|^2 + 2c^2 \left( \cosh(r) - ac \sinh(r) \right) u_t(t,1) d(t) \\
&\quad + \left( 2(1+b+\mu) + M^2 \right) \frac{\cosh(r)}{cr} \|f[t]\|_2^2 - \frac{bkM}{2\lambda} \int_0^1 \theta_x^2(t,x) dx
\end{aligned} \qquad (4.32)$$

Wirtinger's inequality and (2.9) imply that $\int_0^1 \theta_x^2(t,x) dx \geq \pi^2 \int_0^1 \theta^2(t,x) dx$. The previous inequality and the inequality

$$\left( \cosh(r) - ac \sinh(r) \right) u_t(t,1) d(t) \leq \frac{aM}{8} u_t^2(t,1) + \frac{2 \left( \cosh(r) - ac \sinh(r) \right)^2}{aM} |d(t)|^2$$

in conjunction with (4.32) give the following estimate for $t > 0$:

$$\begin{aligned}
\frac{d}{dt} V(u[t], u_t[t], \theta[t]) &\leq -\frac{cr}{8} \int_0^1 \exp(rx) \left( u_t(t,x) + cu_x(t,x) \right)^2 dx \\
&\quad -\frac{cr}{8} \int_0^1 \exp(-rx) \left( u_t(t,x) - cu_x(t,x) \right)^2 dx \\
&\quad -\frac{c}{2} \left( \frac{acM}{2} + \exp(-r)(1+ac)^2 - \exp(r)(1-ac)^2 \right) u_t^2(t,1) \\
&\quad + c^2 \left( c \sinh(r) + \frac{4 \left( \cosh(r) - ac \sinh(r) \right)^2}{aM} + \frac{M}{2a} \right) |d(t)|^2 \\
&\quad + \left( 2(1+b+\mu) + M^2 \right) \frac{\cosh(r)}{cr} \|f[t]\|_2^2 - \frac{bkM}{2\lambda} \pi^2 \int_0^1 \theta^2(t,x) dx
\end{aligned} \qquad (4.33)$$

The fact that $\dfrac{acM}{2} + \exp(-r)(1+ac)^2 - \exp(r)(1-ac)^2 \geq 0$ (recall (4.29)) and definition (4.2) combined with estimate (4.33) give us the following estimate for $t > 0$:

$$\frac{d}{dt}V(u[t],u_t[t],\theta[t]) \leq -\frac{cr}{4}\Phi(u[t],u_t[t]) - \frac{bkM}{2\lambda}\pi^2\|\theta[t]\|_2^2$$
$$+c^2\left(c\sinh(r) + \frac{4(\cosh(r)-ac\sinh(r))^2}{aM} + \frac{M}{2a}\right)|d(t)|^2 \tag{4.34}$$
$$+\left(2(1+b+\mu)+M^2\right)\frac{\cosh(r)}{cr}\|f[t]\|_2^2$$

Definition (4.30), definition (4.22) and estimate (4.3) imply the following inequality:

$$V(u[t],u_t[t],\theta[t]) \leq \left(1+\frac{M}{2}\exp(r)\right)\Phi(u[t],u_t[t]) + \frac{bM}{2\lambda}\|\theta[t]\|_2^2 \tag{4.35}$$

Combining (4.34) and (4.35), we obtain for $t>0$

$$\frac{d}{dt}V(u[t],u_t[t],\theta[t]) \leq -2\omega V(u[t],u_t[t],\theta[t])$$
$$+c^2\left(c\sinh(r) + \frac{M}{2a} + \frac{4(\cosh(r)-ac\sinh(r))^2}{aM}\right)|d(t)|^2 \tag{4.36}$$
$$+\left(2(1+b+\mu)+M^2\right)\frac{\cosh(r)}{cr}\|f[t]\|_2^2$$

where $\omega := \frac{1}{2}\min\left(\frac{cr}{2(2+M\exp(r))}, k\pi^2\right)$. Differential inequality (4.36) directly implies the following estimate for all $t \geq t_0 > 0$:

$$V(u[t],u_t[t],\theta[t]) \leq \exp(-2\omega(t-t_0))V(u[t_0],u_t[t_0],\theta[t_0]) + K_1\sup_{t_0\leq s\leq t}\left(\|f[s]\|_2^2\right) + K_2\sup_{t_0\leq s\leq t}\left(|d(s)|^2\right) \tag{4.37}$$

with $K_1 := \left(2(1+b+\mu)+M^2\right)\frac{\cosh(r)}{2\omega cr}$ and $K_2 := \frac{c^2}{2\omega}\left(\frac{M}{2a} + \frac{4(\cosh(r)-ac\sinh(r))^2}{aM} + c\sinh(r)\right)$.

Since $u \in C^1(\Re_+ \times [0,1])$, $\theta \in C^0(\Re_+ \times [0,1])$ we obtain that $\lim_{t_0\to 0^+}\left(V(u[t_0],u_t[t_0],\theta[t_0])\right) = V(u[0],u_t[0],\theta[0])$ and consequently, we get the following estimate for all $t \geq 0$:

$$V(u[t],u_t[t],\theta[t]) \leq \exp(-2\omega t)V(u[0],u_t[0],\theta[0]) + K_1\sup_{0\leq s\leq t}\left(\|f[s]\|_2^2\right) + K_2\sup_{0\leq s\leq t}\left(|d(s)|^2\right) \tag{4.38}$$

Definition (4.22), definition (4.30) and estimate (4.3) in conjunction with (4.38) allow us to conclude that the following estimate holds for all $t \geq 0$:

$$\min\left(\min\left(1,c^{2}\right)\left(\frac{M}{2}+\exp(-r)\right),\frac{bM}{2\lambda}\right)\left(\|u_{t}[t]\|_{2}^{2}+\|u_{x}[t]\|_{2}^{2}+\|\theta[t]\|_{2}^{2}\right)$$

$$\leq \exp(-2\omega t)\max\left(\max\left(1,c^{2}\right)\left(\frac{M}{2}+\exp(r)\right),\frac{bM}{2\lambda}\right)\left(\|u_{t}[0]\|_{2}^{2}+\|u_{x}[0]\|_{2}^{2}+\|\theta[0]\|_{2}^{2}\right) \quad (4.39)$$

$$+K_{1}\sup_{0\leq s\leq t}\left(\|f[s]\|_{2}^{2}\right)+K_{2}\sup_{0\leq s\leq t}\left(|d(s)|^{2}\right)$$

Estimate (3.2) with appropriate constants $G, \gamma_1, \gamma_2 > 0$, is a direct consequence of estimate (4.39). The proof is complete. ◁

We end this section with the proof of Theorem 3.

**Proof of Theorem 3:** The proof of Theorem 3 uses a similar notation with the proofs of Theorem 1 and Theorem 2 (for example, there is a constant $M > 0$, there are functionals $E, V$, etc.). However, the reader should not be tempted, by an overlapping notation for different quantities, to compare the different quantities in the proofs of Theorem 1, Theorem 2 and Theorem 3. To understand this, we point out that while the functional $E$ in the proof of Theorem 1 is the mechanical energy of the string, the functional $E$ in the proof of Theorem 3 below is not the mechanical energy of the string and is not a functional that may be considered as some kind of energy of the string (e.g., the sum of mechanical energy and the thermal energy of the string).

Consider an arbitrary solution $u \in C^2\left(\Re_+ \times [0,1]\right)$, $\theta \in C^0\left(\Re_+ \times [0,1]\right) \cap C^1\left((0,+\infty); L^2(0,1)\right)$ with $u_{xx} \in C^1\left((0,+\infty); L^2(0,1)\right)$, $\theta[t], u_t[t] \in H^2(0,1)$ for all $t > 0$ of (2.2), (2.8), (2.9), (2.10), (2.11) corresponding to (arbitrary) input $f \in C^0\left(\Re_+; L^2(0,1)\right) \cap C^0\left((0,+\infty) \times (0,1)\right)$. Let $r > 0$ be a constant and define the functionals $E: H^1(0,1) \times C^0([0,1]) \times L^2(0,1) \to \Re_+$, $\Phi: H^1(0,1) \times C^0([0,1]) \to \Re_+$, $W: H^2(0,1) \times L^2(0,1) \to \Re_+$ by means of the formulae

$$E(u,w,\theta) := \frac{1}{2}\int_0^1 w^2(x)dx + \frac{c^2}{2}\int_0^1 u_x^2(x)dx + \frac{b}{2\lambda}\int_0^1 \theta^2(x)dx + \frac{a\sigma}{2}w^2(1),$$

for all $u \in H^2(0,1)$, $w \in C^0([0,1])$, $\theta \in L^2(0,1)$ \hfill (4.40)

$$W(u,w) := \frac{1}{2}\int_0^1 \left(w(x) - \sigma u_{xx}(x)\right)^2 dx, \text{ for all } u \in H^2(0,1), w \in C^0([0,1]) \quad (4.41)$$

$$\Phi(u,w) := \frac{1}{2}\int_0^1 \exp(rx)\left(w(x)+cu_x(x)\right)^2 dx + \frac{1}{2}\int_0^1 \exp(-rx)\left(w(x)-cu_x(x)\right)^2 dx + \frac{B}{2}w^2(1),$$

for all $u \in H^2(0,1)$, $w \in C^0([0,1])$ \hfill (4.42)

where $B := a\sigma\left(\exp(r)(1-ac) + \exp(-r)(1+ac)\right)$. It should be noticed at this point that $E(u[t], u_t[t], \theta[t])$ is an appropriate linear combination at time $t \geq 0$ of the potential energy of the string (the term $\frac{c^2}{2}\int_0^1 u_x^2(t,x)dx$), a measure of the kinetic energy of the string (the term $\frac{1}{2}\int_0^1 u_t^2(t,x)dx + \frac{a\sigma}{2}u_t^2(t,1)$) and the squared $L^2$ norm of the deviation of temperature from its reference value (the term $\int_0^1 \theta^2(t,x)dx$ is not the thermal energy of the string but is the value of a Lyapunov functional used frequently for the heat equation). The functional $W(u,w)$ is a functional that allows the derivation of bounds for the $L^2$ norm of the second spatial derivative of $u$, i.e., $\|u_{xx}\|_2$. The functional $W(u,w)$ was not used in the proofs of Theorem 1 and Theorem 2 because the ISS estimates (3.1), (3.2) do not provide bounds for the $L^2$ norm of the second spatial derivative of $u$. Finally, the functional $\Phi(u,w)$ is not the same functional used in the proofs of Theorem 1 and Theorem 2 (due to the additional term $Bw^2(1)/2$) but when $\sigma \to 0^+$ it becomes equal to the Lyapunov functional for the simple model (A) that was employed in the proofs of Theorem 1 and Theorem 2.

Since $u \in C^2(\Re_+ \times [0,1])$, $\theta \in C^1((0,+\infty); L^2(0,1))$, $u_{xx} \in C^1((0,+\infty); L^2(0,1))$, we get from definitions (4.40), (4.41), (4.42) for all $t > 0$ the following formulae:

$$\frac{d}{dt}E(u[t],u_t[t],\theta[t]) = \int_0^1 u_t(t,x)u_{tt}(t,x)dx + c^2 \int_0^1 u_x(t,x)u_{xt}(t,x)dx \\ + \frac{b}{\lambda}\int_0^1 \theta(t,x)\theta_t(t,x)dx + a\sigma u_t(t,1)u_{tt}(t,1) \quad (4.43)$$

$$\frac{d}{dt}W(u[t],u_t[t]) = \int_0^1 \left(u_t(t,x) - \sigma u_{xx}(t,x)\right)\left(u_{tt}(t,x) - \sigma u_{xxt}(t,x)\right)dx \quad (4.44)$$

$$\frac{d}{dt}\Phi(u[t],u_t[t]) = \int_0^1 \exp(rx)\left(u_t(t,x) + cu_x(t,x)\right)\left(u_{tt}(t,x) + cu_{xt}(t,x)\right)dx \\ + \int_0^1 \exp(-rx)\left(u_t(t,x) - cu_x(t,x)\right)\left(u_{tt}(t,x) - cu_{xt}(t,x)\right)dx + Bu_t(t,1)u_{tt}(t,1) \quad (4.45)$$

Boundary condition (2.11) and the fact that $u \in C^2(\Re_+ \times [0,1])$ implies that

$$u_{xt}(t,1) = -au_{tt}(t,1), \text{ for } t > 0. \quad (4.46)$$

Integrating by parts the integral $\int_0^1 u_x(t,x)u_{xt}(t,x)dx$ and using (4.43), (2.8), (2.10), (2.11), (2.2) (the latter implies that $u_t(t,0) = 0$ for all $t > 0$) we obtain for $t > 0$:

$$\frac{d}{dt}E(u[t],u_t[t],\theta[t]) = -\mu \int_0^1 u_t^2(t,x)dx + \int_0^1 u_t(t,x)f(t,x)dx + \sigma \int_0^1 u_t(t,x)u_{xxt}(t,x)dx \\ -ac^2u_t^2(t,1) - b\int_0^1 \left(u_t(t,x)\theta(t,x)\right)_x dx + \frac{kb}{\lambda}\int_0^1 \theta(t,x)\theta_{xx}(t,x)dx + a\sigma u_t(t,1)u_{tt}(t,1) \quad (4.47)$$

Integrating by parts the integrals $\int_0^1 u_t(t,x)u_{xxt}(t,x)dx$, $\int_0^1 \theta(t,x)\theta_{xx}(t,x)dx$ and using (4.46), (4.47), (2.9), (2.2) (the latter implies that $u_t(t,0) = 0$ for all $t > 0$) we obtain for $t > 0$:

$$\frac{d}{dt}E(u[t],u_t[t],\theta[t]) = -\mu\int_0^1 u_t^2(t,x)dx + \int_0^1 u_t(t,x)f(t,x)dx$$
$$-\sigma\int_0^1 u_{xt}^2(t,x)dx - ac^2 u_t^2(t,1) - \frac{kb}{\lambda}\int_0^1 \theta_x^2(t,x)dx \qquad (4.48)$$

Using (4.45) and (2.10) we get for $t > 0$:

$$\frac{d}{dt}\Phi(u[t],u_t[t]) = c\int_0^1 \exp(rx)\left(\frac{1}{2}(u_t(t,x)+cu_x(t,x))^2\right)_x dx$$
$$-c\int_0^1 \exp(-rx)\left(\frac{1}{2}(u_t(t,x)-cu_x(t,x))^2\right)_x dx + Bu_t(t,1)u_{tt}(t,1)$$
$$+\int_0^1 \exp(rx)(u_t(t,x)+cu_x(t,x))(f(t,x)+\sigma u_{xxt}(t,x)-b\theta_x(t,x)-\mu u_t(t,x))dx \qquad (4.49)$$
$$+\int_0^1 \exp(-rx)(u_t(t,x)-cu_x(t,x))(f(t,x)+\sigma u_{xxt}(t,x)-b\theta_x(t,x)-\mu u_t(t,x))dx$$

Integrating by parts the integrals $\int_0^1 \exp(rx)\left(\frac{1}{2}(u_t(t,x)+cu_x(t,x))^2\right)_x dx$, $\int_0^1 \exp(-rx)\left(\frac{1}{2}(u_t(t,x)-cu_x(t,x))^2\right)_x dx$ and using (4.49), (2.11), (2.2) (the latter implies that $u_t(t,0) = 0$ for all $t > 0$) we obtain for $t > 0$:

$$\frac{d}{dt}\Phi(u[t],u_t[t]) = -\frac{cr}{2}\int_0^1 \exp(rx)(u_t(t,x)+cu_x(t,x))^2 dx$$
$$-\frac{cr}{2}\int_0^1 \exp(-rx)(u_t(t,x)-cu_x(t,x))^2 dx$$
$$+\frac{c}{2}\left(\exp(r)(1-ac)^2 - \exp(-r)(1+ac)^2\right)u_t^2(t,1) + Bu_t(t,1)u_{tt}(t,1) \qquad (4.50)$$
$$+\int_0^1 \exp(rx)(u_t(t,x)+cu_x(t,x))(f(t,x)+\sigma u_{xxt}(t,x)-b\theta_x(t,x)-\mu u_t(t,x))dx$$
$$+\int_0^1 \exp(-rx)(u_t(t,x)-cu_x(t,x))(f(t,x)+\sigma u_{xxt}(t,x)-b\theta_x(t,x)-\mu u_t(t,x))dx$$

Integrating by parts the integrals $\int_0^1 \exp(\pm rx)(u_t(t,x)\pm cu_x(t,x))u_{xxt}(t,x)dx$ and using (4.46), (4.50), (2.11), (2.2) (the latter implies that $u_t(t,0) = 0$ for all $t > 0$) and the fact that $B = a\sigma\left(\exp(r)(1-ac)+\exp(-r)(1+ac)\right)$, we obtain for $t > 0$:

$$\frac{d}{dt}\Phi(u[t],u_t[t]) = -\frac{cr}{2}\int_0^1 \exp(rx)\big(u_t(t,x)+cu_x(t,x)\big)^2 dx$$

$$-\frac{cr}{2}\int_0^1 \exp(-rx)\big(u_t(t,x)-cu_x(t,x)\big)^2 dx$$

$$-\frac{c}{2}\Big(\exp(-r)(1+ac)^2 - \exp(r)(1-ac)^2\Big)u_t^2(t,1)$$

$$-2\sigma\int_0^1 \cosh(rx)u_{xt}^2(t,x)dx - 2c\sigma\int_0^1 \sinh(rx)u_{xx}(t,x)u_{xt}(t,x)dx \qquad (4.51)$$

$$+\int_0^1 \exp(rx)\big(u_t(t,x)+cu_x(t,x)\big)\big(f(t,x)-\sigma r u_{xt}(t,x)-b\theta_x(t,x)-\mu u_t(t,x)\big)dx$$

$$+\int_0^1 \exp(-rx)\big(u_t(t,x)-cu_x(t,x)\big)\big(f(t,x)+\sigma r u_{xt}(t,x)-b\theta_x(t,x)-\mu u_t(t,x)\big)dx$$

Using the following inequalities

$$\big(u_t(t,x)\pm cu_x(t,x)\big)f(t,x) \leq \frac{cr}{4(1+\mu+b+\sigma r)}\big(u_t(t,x)\pm cu_x(t,x)\big)^2 + \frac{1+\mu+b+\sigma r}{cr}|f(t,x)|^2$$

$$\big(u_t(t,x)\pm cu_x(t,x)\big)\theta_x(t,x) \geq -\frac{cr}{4(1+\mu+b+\sigma r)}\big(u_t(t,x)\pm cu_x(t,x)\big)^2 - \frac{1+\mu+b+\sigma r}{cr}\theta_x^2(t,x)$$

$$\big(u_t(t,x)\pm cu_x(t,x)\big)u_t(t,x) \geq -\frac{cr}{4(1+\mu+b+\sigma r)}\big(u_t(t,x)\pm cu_x(t,x)\big)^2 - \frac{1+\mu+b+\sigma r}{cr}u_t^2(t,x)$$

$$\big|u_t(t,x)\pm cu_x(t,x)\big|\big|u_{xt}(t,x)\big| \leq \frac{cr}{4(1+\mu+b+\sigma r)}\big(u_t(t,x)\pm cu_x(t,x)\big)^2 + \frac{1+\mu+b+\sigma r}{cr}u_{xt}^2(t,x)$$

we obtain from (4.51) the following estimate for $t>0$:

$$\frac{d}{dt}\Phi(u[t],u_t[t]) \leq -\frac{cr}{4}\int_0^1 \exp(rx)\big(u_t(t,x)+cu_x(t,x)\big)^2 dx$$

$$-\frac{cr}{4}\int_0^1 \exp(-rx)\big(u_t(t,x)-cu_x(t,x)\big)^2 dx$$

$$-\frac{c}{2}\Big(\exp(-r)(1+ac)^2 - \exp(r)(1-ac)^2\Big)u_t^2(t,1) \qquad (4.52)$$

$$-2c\sigma\int_0^1 \sinh(rx)u_{xx}(t,x)u_{xt}(t,x)dx$$

$$+\frac{2(1+\mu+b+\sigma r)}{cr}\int_0^1 \cosh(rx)\Big(|f(t,x)|^2 + b\theta_x^2(t,x) + \mu u_t^2(t,x) + \sigma r u_{xt}^2(t,x)\Big)dx$$

Equation (4.44) in conjunction with (2.10) implies for all $t>0$:

$$\frac{d}{dt}W(u[t],u_t[t]) = \int_0^1 \big(u_t(t,x) - \sigma u_{xx}(t,x)\big)\big(c^2 u_{xx}(t,x) - \mu u_t(t,x) - b\theta_x(t,x) + f(t,x)\big)dx$$

$$= -\sigma c^2 \int_0^1 u_{xx}^2(t,x)dx + (c^2 + \mu\sigma)\int_0^1 u_t(t,x)u_{xx}(t,x)dx - \mu\int_0^1 u_t^2(t,x)dx \qquad (4.53)$$

$$-b\int_0^1 u_t(t,x)\theta_x(t,x)dx + \int_0^1 u_t(t,x)f(t,x)dx + \sigma b\int_0^1 u_{xx}(t,x)\theta_x(t,x)dx - \sigma\int_0^1 u_{xx}(t,x)f(t,x)dx$$

Using the inequalities

$$u_t(t,x)u_{xx}(t,x) \le \frac{\sigma c^2}{2(c^2 + (1+\mu+b)\sigma)}u_{xx}^2(t,x) + \frac{c^2 + (1+\mu+b)\sigma}{2\sigma c^2}u_t^2(t,x)$$

$$u_{xx}(t,x)\theta_x(t,x) \le \frac{\sigma c^2}{2(c^2 + (1+\mu+b)\sigma)}u_{xx}^2(t,x) + \frac{c^2 + (1+\mu+b)\sigma}{2\sigma c^2}\theta_x^2(t,x)$$

$$u_{xx}(t,x)f(t,x) \ge -\frac{\sigma c^2}{2(c^2 + (1+\mu+b)\sigma)}u_{xx}^2(t,x) - \frac{c^2 + (1+\mu+b)\sigma}{2\sigma c^2}|f(t,x)|^2$$

$$u_t(t,x)\theta_x(t,x) \ge -\frac{\sigma c^2}{2(c^2 + (1+\mu+b)\sigma)}u_t^2(t,x) - \frac{c^2 + (1+\mu+b)\sigma}{2\sigma c^2}\theta_x^2(t,x)$$

$$u_t(t,x)f(t,x) \le \frac{\sigma c^2}{2(c^2 + (1+\mu+b)\sigma)}u_t^2(t,x) + \frac{c^2 + (1+\mu+b)\sigma}{2\sigma c^2}|f(t,x)|^2$$

we obtain from (4.53) the following estimate for $t > 0$:

$$\frac{d}{dt}W(u[t],u_t[t]) \le -\frac{\sigma c^2}{2}\int_0^1 u_{xx}^2(t,x)dx$$

$$+ Q\int_0^1 u_t^2(t,x)dx + \frac{\mu(1+\mu+b)\sigma}{2c^2}\int_0^1 u_t^2(t,x)dx \qquad (4.54)$$

$$+ \frac{(\sigma+1)b}{2\sigma c^2}\big(c^2 + (1+\mu+b)\sigma\big)\int_0^1 \theta_x^2(t,x)dx + \frac{\sigma+1}{2\sigma c^2}\big(c^2 + (1+\mu+b)\sigma\big)\int_0^1 |f(t,x)|^2 dx$$

where $Q := \frac{c^2 + (1+b)\sigma}{2\sigma} + \frac{\sigma c^2(1+b)}{2(c^2 + (1+\mu+b)\sigma)}$. Define:

$$R := \frac{cr}{8Q\cosh(r)} \qquad (4.55)$$

$$V(u,w,\theta) := \Phi(u,w) + RW(u,w) + M\,E(u,w,\theta),$$
$$\text{for all } u \in H^2(0,1),\ w \in C^0([0,1]),\ \theta \in L^2(0,1) \qquad (4.56)$$

where $M > 0$ is a sufficiently large constant that satisfies

$$M \geq \frac{2(1+\mu+b+\sigma r)}{c}\cosh(r) + \frac{4\sinh^2(r)}{cR}$$

$$M \geq \frac{(1+\mu+b)\sigma R}{2c^2} + \frac{2(1+\mu+b+\sigma r)}{cr}\cosh(r)$$

$$M \geq \frac{2\lambda}{k}\left(\frac{2(1+\mu+b+\sigma r)}{cr}\cosh(r) + \frac{(\sigma+1)}{2\sigma c^2}\left(c^2 + (1+\mu+b)\sigma\right)R\right) \quad (4.57)$$

$$M \geq \frac{1}{ac}\left(\exp(r)(1-ac)^2 - \exp(-r)(1+ac)^2\right)$$

$$M > R - 2\exp(-r)$$

$$M > \frac{-B}{a\sigma}$$

Definitions (4.40), (4.41), (4.42) and (4.56) guarantee that the following inequalities hold for all $u \in H^2(0,1)$, $w \in C^0([0,1])$, $\theta \in L^2(0,1)$:

$$V(u,w,\theta) \leq C_2\left(\|w\|_2^2 + \|u_x\|_2^2 + \|\theta\|_2^2 + w^2(1) + \|u_{xx}\|_2^2\right)$$
$$V(u,w,\theta) \geq C_1\left(\|w\|_2^2 + \|u_x\|_2^2 + \|\theta\|_2^2 + w^2(1) + \|u_{xx}\|_2^2\right) \quad (4.58)$$

where

$$C_1 := \min\left(\frac{M-R}{2} + \exp(-r), c^2\left(\frac{M}{2} + \exp(-r)\right), \frac{bM}{2\lambda}, \frac{a\sigma M + B}{2}, \frac{R\sigma^2}{4}\right)$$

$$C_2 := \max\left(\frac{M}{2} + \exp(r) + R, c^2\left(\frac{M}{2} + \exp(r)\right), \frac{bM}{2\lambda}, \frac{a\sigma M + B}{2}, R\sigma^2\right) \quad (4.59)$$

Notice that due to (4.57), the constants $C_1, C_2$ defined by (4.59) are positive, i.e., $0 < C_1 \leq C_2$.

Using (4.48), (4.52), (4.54) and definition (4.56), we get for $t > 0$:

$$\frac{d}{dt}V(u[t], u_t[t], \theta[t]) \leq -\frac{cr}{4}\int_0^1 \exp(rx)\left(u_t(t,x) + cu_x(t,x)\right)^2 dx$$

$$-\frac{cr}{4}\int_0^1 \exp(-rx)\left(u_t(t,x) - cu_x(t,x)\right)^2 dx$$

$$-\frac{c}{2}\left(2acM + \exp(-r)(1+ac)^2 - \exp(r)(1-ac)^2\right)u_t^2(t,1)$$

$$-\left(\sigma M - \frac{2\sigma(1+\mu+b+\sigma r)}{c}\cosh(r)\right)\int_0^1 u_{xt}^2(t,x)dx - \frac{\sigma c^2 R}{2}\int_0^1 u_{xx}^2(t,x)dx$$

$$+ M\int_0^1 u_t(t,x)f(t,x)dx + QR\int_0^1 u_t^2(t,x)dx - 2c\sigma\int_0^1 \sinh(rx)u_{xx}(t,x)u_{xt}(t,x)dx$$

$$+\left(\frac{\sigma+1}{2\sigma c^2}\left(c^2 + (1+\mu+b)\sigma\right)R + \frac{2(1+\mu+b+\sigma r)}{cr}\cosh(r)\right)\|f[t]\|_2^2$$

$$-\mu\left(M - \frac{(1+\mu+b)\sigma R}{2c^2} - \frac{2(1+\mu+b+\sigma r)}{cr}\cosh(r)\right)\int_0^1 u_t^2(t,x)dx$$

$$-b\left(\frac{kM}{\lambda} - \frac{2(1+\mu+b+\sigma r)}{cr}\cosh(r) - \frac{(\sigma+1)}{2\sigma c^2}\left(c^2 + (1+\mu+b)\sigma\right)R\right)\int_0^1 \theta_x^2(t,x)dx$$

$$(4.60)$$

Using inequalities (4.57) and the inequalities

$$Mu_t(t,x)f(t,x) \le QRu_t^2(t,x) + \frac{M^2}{4QR}|f(t,x)|^2$$

$$u_{xx}(t,x)u_{xt}(t,x) \ge -\frac{cR}{8\sinh(r)}u_{xx}^2(t,x) - \frac{2\sinh(r)}{cR}u_{xt}^2(t,x)$$

we obtain from (4.60) for all $t > 0$:

$$\begin{aligned}\frac{d}{dt}V(u[t],u_t[t],\theta[t]) \le &-\frac{cr}{4}\int_0^1 \exp(rx)\bigl(u_t(t,x)+cu_x(t,x)\bigr)^2 dx\\ &-\frac{cr}{4}\int_0^1 \exp(-rx)\bigl(u_t(t,x)-cu_x(t,x)\bigr)^2 dx - \frac{ac^2M}{2}u_t^2(t,1)\\ &-\sigma\left(M - \frac{2(1+\mu+b+\sigma r)}{c}\cosh(r) - \frac{4\sinh^2(r)}{cR}\right)\int_0^1 u_{xt}^2(t,x)dx\\ &-\frac{\sigma c^2 R}{4}\int_0^1 u_{xx}^2(t,x)dx + 2QR\int_0^1 u_t^2(t,x)dx - \frac{bkM}{2\lambda}\int_0^1 \theta_x^2(t,x)dx + K\|f[t]\|_2^2\end{aligned} \quad (4.61)$$

where $K := \frac{\sigma+1}{2\sigma c^2}\bigl(c^2+(1+\mu+b)\sigma\bigr)R + \frac{2(1+\mu+b+\sigma r)}{cr}\cosh(r) + \frac{M^2}{4QR}$. Combining (4.11) with (4.61) and using definition (4.55) and inequalities (4.57), we obtain the following estimate for $t > 0$:

$$\begin{aligned}\frac{d}{dt}V(u[t],u_t[t],\theta[t]) \le &-\frac{cr}{8}\int_0^1 \exp(rx)\bigl(u_t(t,x)+cu_x(t,x)\bigr)^2 dx\\ &-\frac{cr}{8}\int_0^1 \exp(-rx)\bigl(u_t(t,x)-cu_x(t,x)\bigr)^2 dx - \frac{ac^2M}{2}u_t^2(t,1)\\ &-\frac{\sigma c^2 R}{4}\|u_{xx}[t]\|_2^2 - \frac{bkM}{2\lambda}\|\theta_x[t]\|_2^2 + K\|f[t]\|_2^2\end{aligned} \quad (4.62)$$

Wirtinger's inequality and (2.9) imply that $\int_0^1 \theta_x^2(t,x)dx \ge \pi^2 \int_0^1 \theta^2(t,x)dx$. The previous inequality and estimate (4.62) give the following estimate for $t > 0$:

$$\begin{aligned}\frac{d}{dt}V(u[t],u_t[t],\theta[t]) \le &-\frac{cr}{8}\exp(-r)\|u_t[t]\|_2^2 - \frac{c^3 r}{8}\exp(-r)\|u_x[t]\|_2^2\\ &-\frac{ac^2M}{2}u_t^2(t,1) - \frac{\sigma c^2 R}{4}\|u_{xx}[t]\|_2^2 - \frac{bkM}{2\lambda}\pi^2\|\theta[t]\|_2^2 + K\|f[t]\|_2^2\end{aligned} \quad (4.63)$$

Consequently, we obtain from (4.63) for $t > 0$:

$$\frac{d}{dt}V(u[t],u_t[t],\theta[t]) \le -\varphi\Bigl(\|u_t[t]\|_2^2 + \|u_x[t]\|_2^2 + u_t^2(t,1) + \|u_{xx}[t]\|_2^2 + \|\theta[t]\|_2^2\Bigr) + K\|f[t]\|_2^2 \quad (4.64)$$

where $\varphi := \min\left(\frac{cr}{8}\exp(-r), \frac{c^3 r}{8}\exp(-r), \frac{ac^2 M}{2}, \frac{\sigma c^2 R}{4}, \frac{bkM}{2\lambda}\pi^2\right)$. Using (4.58) and (4.64) we get for $t > 0$:

$$\frac{d}{dt}V(u[t],u_t[t],\theta[t]) \leq -2\omega V(u[t],u_t[t],\theta[t]) + K\|f[t]\|_2^2 \qquad (4.65)$$

where $\omega := \frac{\varphi}{2C_2}$. Differential inequality (4.65) directly implies the following estimate for all $t \geq t_0 > 0$:

$$V(u[t],u_t[t],\theta[t]) \leq \exp(-2\omega(t-t_0))V(u[t_0],u_t[t_0],\theta[t_0]) + \frac{K}{2\omega}\sup_{t_0 \leq s \leq t}\left(\|f[s]\|_2^2\right) \qquad (4.66)$$

Since $u \in C^2(\Re_+ \times [0,1])$, $\theta \in C^0(\Re_+ \times [0,1])$ we obtain that $\lim_{t_0 \to 0^+}\left(V(u[t_0],u_t[t_0],\theta[t_0])\right) = V(u[0],u_t[0],\theta[0])$ and consequently, we get the following estimate for all $t \geq 0$:

$$V(u[t],u_t[t],\theta[t]) \leq \exp(-2\omega t)V(u[0],u_t[0],\theta[0]) + \frac{K}{2\omega}\sup_{0 \leq s \leq t}\left(\|f[s]\|_2^2\right) \qquad (4.67)$$

Estimate (3.3) with appropriate constants $G, \gamma > 0$, is a direct consequence of estimate (4.67) and inequalities (4.58). The proof is complete. ◁

**Remark:** It should be noticed that when $\sigma \to 0^+$, the constants $K := \frac{\sigma+1}{2\sigma c^2}\left(c^2 + (1+\mu+b)\sigma\right)R + \frac{2(1+\mu+b+\sigma r)}{cr}\cosh(r) + \frac{M^2}{4QR}$, $\omega := \frac{\varphi}{2C_2}$ and $C_1$ defined by (4.59) satisfy $C_1 \to 0^+$, $K \to +\infty$ and $\omega \to 0^+$. Therefore, it becomes clear from (4.67) and (4.58) that the gain $\gamma$ in the ISS estimate (3.3) of the distributed disturbance input $f$ becomes unbounded as $\sigma \to 0^+$. This is expected because when $\sigma \to 0^+$, model (D) "tends" to model (C) for which the ISS estimate (3.2) does not allow the derivation of bounds for $\|u_{xx}[t]\|_2$ and $|u_t(t,1)|$. On the other hand, the asymptotic analysis when $\sigma \to +\infty$ of the gain $\gamma$ in the ISS estimate (3.3) of the input $f$ is not easy and is a topic for future research.

## 5. Concluding Remarks

The study of the robustness properties of the 1-D wave equation for an elastic vibrating string was performed under four different damping mechanisms that are usually neglected in the study of the wave equation: (i) friction with the surrounding medium of the string (or viscous damping), (ii) thermoelastic phenomena (or thermal damping), (iii) internal friction of the string (or Kelvin-Voigt damping), and (iv) friction at the free end of the string (the so-called passive damper).

The study is by no means complete. Future work may consider different boundary conditions for the temperature at the ends of the string or may also consider the heat exchange between the body of the string and the surrounding fluid. Moreover, from a control perspective, the "ultimate goal" would be to use the ISS analysis for the design of robust stabilizers (robust even in the presence of small input delays). We cannot be sure that such a goal is feasible (although some stabilization results are given in [24]) but it is certainly a topic that must be studied in the future.